\documentclass{article}%
\usepackage{amssymb}
\usepackage{amsfonts}
\usepackage{amsmath}
\usepackage{graphicx}%
\newtheorem{theorem}{Theorem}

\newtheorem{example}[theorem]{Example}

\newtheorem{lemma}[theorem]{Lemma}

\newtheorem{proposition}[theorem]{Proposition}
\newtheorem{remark}[theorem]{Remark}

\newenvironment{proof}[1][Proof]{\noindent\textbf{#1.} }{\ \rule{0.5em}{0.5em}}
\begin{document}

\title{$ADE$ bundles over $ADE$ singular surfaces and  \\flag varieties of $ADE$ type}
\author{Yunxia Chen \& Naichung Conan Leung}
\date{}
\maketitle

\begin{abstract}
Based on the Brieskorn-Slodowy-Grothendieck diagram, we write the holomorphic structures (or filtrations) of the ADE Lie algebra bundles over the corresponding type ADE flag varieties, over the cotangent bundles of these flag varieties, and over the corresponding type $ADE$ singular surfaces. The main tool is the cohomology of line bundles over flag varieties and their cotangent bundles.
\end{abstract}

\section{Introduction}
This paper is a continuation of our earlier paper about $ADE$ bundles over $ADE$ singular surfaces \cite{CL}. In that paper, for every $ADE$ singular compact surface with $p_{g}=0$, we constructed a corresponding type $ADE$ bundle over it, using the exceptional locus in its minimal resolution and bundle extensions. There are lots of studies for bundles over surfaces (\cite{C}\cite{Don}\cite{Lee}\cite{LXZ}\cite{LZ2}\cite{Loo}). In this paper, we will base on Slodowy's paper \cite{Slo}, study the homogenous $ADE$ bundles over flag varieties of the corresponding type, and their lifts to the cotangent bundles of the flag varieties, and then their restrictions to the $ADE$ singular surfaces of the corresponding type.

In more detail, given a complex simple Lie algebra $\mathfrak{g}$ of $ADE$ type (we will use $\mathfrak{b}$, $\mathfrak{n}$, $\mathfrak{t}$, $G$, $B$, $W$ to denote the corresponding standard lower-triangular Borel subalgebra, standard lower-triangular nilpotent subalgebra, standard Cartan subalgebra, simply-connected Lie group, standard lower-triangular Borel subgroup, Weyl group respectively), we can have an $ADE$ singular non-compact surface $S$ of the corresponding type, as the intersection of the transversal slice $S_{x}$ of a subregular nilpotent element $x$ and the nilpotent variety $N(\mathfrak{g})$ of $\mathfrak{g}$. Furthermore, the restriction of the adjoint quotient $\mathfrak{g}\rightarrow \mathfrak{t}/W$ to the transversal slice $S_{x}$ is a semiuniversal deformation of the corresponding $ADE$ singularity. This result is conjectured by Grothendieck and proved by Brieskorn in 1970 \cite{Bri}. After that, Grothendieck defined a morphism $G \times \mathfrak{b}/B\rightarrow \mathfrak{t}$ and gave a simultaneous resolution of the adjoint quotient $\mathfrak{g}\rightarrow \mathfrak{t}/W$ using it. The restriction of the Grothendieck resolution to the above transversal slice $S_{x}$ is also a simultaneous resolution \cite{Slo}. In 1969, Springer gave a resolution of singularities for the nilpotent variety $N(\mathfrak{g})$ through $G \times \mathfrak{n}/B\rightarrow N(\mathfrak{g})$, note that $G \times \mathfrak{n}/B \cong T^{*}(G/B)$ is the cotangent bundle of the flag variety $G/B$. The connection among these resolutions can be shown in the following Brieskorn-Slodowy-Grothendieck diagram (here $\widetilde{S}$ is the minimal resolution of $S$ and $C=\bigcup C_i$ is the exceptional locus with each $C_i$ irreducible component).

\[
\begin{array}{l}
 \begin{array}{*{20}c}
  {C=\bigcup C_i} & {\subset} & {\widetilde{S}} & {\longrightarrow} & {S= N(\mathfrak{g}) \cap S_x }  \\
    \cap  & {} & \cap  & {}  &  \cap   \\
 {G/B} & {\subset} &  {G \times \mathfrak{n}/B} & {\longrightarrow} & {N(\mathfrak{g})}  \\
   {} & {} & \cap  & {\kern 1pt}  &  \cap   \\
 {} & {} & {G \times \mathfrak{b}/B} & {\longrightarrow} & {\mathfrak{g}}  \\
    {} & {} &  \downarrow  & {\kern 1pt}  &  \downarrow  \\
     {} & {} &  {\mathfrak{t}} & {\longrightarrow} & {\mathfrak{t}/W}  \\
\end{array} \\
\end{array}
\]

\bigskip
Given the above background, we want to understand the associated Lie algebra bundles $G \times \mathfrak{g}/B$ over $G/B$, and $G \times \mathfrak{n}\times \mathfrak{g}/B$ over $T^*(G/B)$ respectively. It is obvious that these bundles are trivial as the action of $B$ on $\mathfrak{g}$ can extend to the whole $G$.
What we want to do is to describe natural holomorphic filtration structures on these bundles explicitly. Since the minimal resolution $\widetilde{S}$ of the $ADE$ singular surface $S$ is contained in $G \times \mathfrak{n}/B$, we can consider the restriction of the $\mathfrak{g}$-bundle $G \times \mathfrak{n}\times \mathfrak{g}/B$ from $G \times \mathfrak{n}/B$ to $\widetilde{S}$. We will rewrite this $\mathfrak{g}$-bundle over $\widetilde{S}$ in terms of the exceptional locus of $\widetilde{S}$, and compare it with the $ADE$ bundle we constructed in \cite{CL}.

The organization of this paper is as follows. Section 2 gives a quick review of the construction of $ADE$ bundles over $ADE$ singular surfaces \cite{CL}. In section 3, we describe the filtrations of $G \times \mathfrak{g}/B$ over $G/B$. In section 4, we describe the filtrations of $G \times \mathfrak{n}\times \mathfrak{g}/B$ over $T^*(G/B)$. Section 5 describes the restriction of $G \times \mathfrak{n}\times \mathfrak{g}/B$ to the minimal resolution of the $ADE$ singular surface.

\bigskip

$\mathbf{Acknowledgements.}$ We are grateful to
J.J. Zhang for many useful comments and discussions.
The first author is supported by the National Natural
Science Foundation of China (No. 11501201), the China Postdoctoral Science Foundation (No. 2015M570334) and the Fundamental Research Funds for the
Central Universities under project no. 222201514322.
The second
author is supported by a research grant from the Research Grants Council of
the Hong Kong Special Administrative Region, China (No. 2130405).

\section{$ADE$ bundles over $ADE$ singular surfaces}
In this section, we review the construction of $ADE$ bundles over $ADE$ singular surfaces with $p_{g}=0$.

An $ADE$ singularity in a surface $X$ can be described locally as a
quotient singularity $\mathbb{C}^{2}/\Gamma$ with $\Gamma$ a finite subgroup
of $SL(2,\mathbb{C)}$. It is also called a Kleinian singularity or simple
singularity \cite{BPV}. If
we consider the minimal resolution $Y$ of $X$, then every irreducible
component of the exceptional locus $C=\bigcup C_{i}$ is a smooth
rational curve with normal bundle $\mathcal{O}_{\mathbb{P}^{1}}(-2)$, i.e. a
$(-2)$-curve, and the dual graph of the exceptional locus is an $ADE$ Dynkin
diagram.

There is a natural decomposition%
\[
H^{2}(Y,\mathbb{Z})=H^{2}(X,\mathbb{Z})\oplus\Pi\text{,}%
\]
where $\Pi=\{\sum a_{i}[C_{i}]|a_{i}\in\mathbb{Z}\}$. The set
$\Phi:=\{\alpha\in\Pi|\alpha^{2}=-2\}$ is a simply-laced (i.e. $\ ADE$)
root system of a simple Lie algebra $\mathfrak{g}$ and $\Delta=\{[C_{i}]\}$ is
a base of $\Phi$. For any $\alpha\in\Phi$, there exists a unique divisor
$D=\sum a_{i}C_{i}$ with $\alpha=[D]$, and we define a line bundle
$\mathcal{O}(\alpha):=\mathcal{O}(D)$ over $Y$.

We define a Lie algebra bundle of type $\mathfrak{g}$ over $Y$ as follows:

\begin{center}
$\mathcal{E}_{0}^{\mathfrak{g}}:=\mathcal{O}^{\oplus r}\oplus\bigoplus_{\alpha\in\Phi
}\mathcal{O}(\alpha)$.
\end{center}

For every open chart $U$ of $Y$, we take $x_{\alpha}^{U}$ to be a nonvanishing
holomorphic section of $\mathcal{O}_{U}(\alpha)$ and $h_{i}^{U}$ ($i=1,\cdots,r$)
nonvanishing holomorphic sections of $\mathcal{O}_{U}^{\oplus r}$. Define a Lie algebra
structure $[ ~ , ~ ]$ on $\mathcal{E}_{0}^{\mathfrak{g}}$ such that $\{x_{\alpha
}^{U}$'s, $h_{i}^{U}$'s$\}$ is the Chevalley basis \cite{H2}, i.e.

(a) $[h_{i}^{U},$ $h_{j}^{U}]=0$, $1\leq i$, $j\leq r$.

(b) $[h_{i}^{U},$ $x_{\alpha}^{U}]=\left\langle \alpha\text{, }C_{i}%
\right\rangle x_{\alpha}^{U}$, $1\leq i\leq r$, $\alpha\in\Phi$.

(c) $[x_{\alpha}^{U},$ $x_{-\alpha}^{U}]=h_{\alpha}^{U}$ is a $\mathbb{Z}%
$-linear combination of $h_{i}^{U}$.

(d) If $\alpha$, $\beta$ are independent roots, and $\beta-p\alpha
,\cdots,\beta+q\alpha$ is the $\alpha$-string through $\beta$, then
$[x_{\alpha}^{U},$ $x_{\beta}^{U}]=0$ if $q=0$, otherwise $[x_{\alpha}^{U},$
$x_{\beta}^{U}]=\pm(p+1)x_{\alpha+\beta}^{U}$.

Since $\mathfrak{g}$ is simply-laced, all its roots have the same length, we
have any $\alpha$-string through $\beta$ is of length at most $2$. So (d) can
be written as $[x_{\alpha}^{U},$ $x_{\beta}^{U}]=n_{\alpha,\beta}%
x_{\alpha+\beta}^{U}$, where $n_{\alpha,\beta}=\pm1$ if $\alpha+\beta\in\Phi$,
otherwise $n_{\alpha,\beta}=0$. From the Jacobi identity, we have for any
$\alpha,\beta,\gamma\in\Phi$, $n_{\alpha,\beta}n_{\alpha+\beta,\gamma
}+n_{\beta,\gamma}n_{\beta+\gamma,\alpha}+n_{\gamma,\alpha}n_{\gamma
+\alpha,\beta}=0$. This Lie algebra structure is compatible with different
trivializations of $\mathcal{E}_{0}^{\mathfrak{g}}$ \cite{LZ}.

By Friedman-Morgan \cite{FM}, a bundle over $Y$ can descend to $X$ if and only
if its restriction to each irreducible component $C_{i}$ of the exceptional
locus is trivial. But $\mathcal{E}_{0}^{\mathfrak{g}}|_{C_{i}}$ is not
trivial as $\mathcal{O}(C_{i})|_{C_{i}}\cong \mathcal{O}_{\mathbb{P}^{1}}(-2)$. We will
construct a new holomorphic structure on $\mathcal{E}_{0}^{\mathfrak{g}}$,
which preserves the Lie algebra structure and the resulting bundle
$\mathcal{E}_{\varphi}^{\mathfrak{g}}$ can descend to $X$.

As we have fixed a base $\Delta$ of $\Phi$, we have a decomposition $\Phi
=\Phi^{+}\cup\Phi^{-}$ into positive and negative roots.

\bigskip
$\mathbf{Definition}$
Given any $\varphi=(\varphi_{\alpha})_{\alpha\in\Phi^{+}}\in\Omega
^{0,1}(Y,\bigoplus_{\alpha\in\Phi^{+}}\mathcal{O}(\alpha))$, we define $\overline
{\partial}_{\varphi}:\Omega^{0,0}(Y,\mathcal{E}_{0}^{\mathfrak{g}%
})\longrightarrow\Omega^{0,1}(Y,\mathcal{E}_{0}^{\mathfrak{g}})$ by
\[
\overline{\partial}_{\varphi}:=\overline{\partial}_{0}+ad(\varphi
):=\overline{\partial}_{0}+\sum_{\alpha\in\Phi^{+}}ad(\varphi_{\alpha
})\text{,}%
\]
where $\overline{\partial}_{0}$ is the standard holomorphic structure of
$\mathcal{E}_{0}^{\mathfrak{g}}$. More explicitly, if we write $\varphi
_{\alpha}=c_{\alpha}^{U}x_{\alpha}^{U}$ locally for some one form $c_{\alpha
}^{U}$, then $ad(\varphi_{\alpha})=c_{\alpha}^{U}ad(x_{\alpha}^{U})$.
\bigskip

From the Jacobi identity, we have $\overline{\partial}_{\varphi}$ is compatible with the Lie
algebra structure, i.e. $\overline{\partial}_{\mathcal{\varphi}}[ ~ , ~ ]=0$.

For $\overline{\partial}_{\mathcal{\varphi}}$ to define a holomorphic
structure, we need
\[
0=\overline{\partial}_{\mathcal{\varphi}}^{2}=\sum_{\alpha\in\Phi^{+}%
}(\overline{\partial}_{0}c_{\alpha}^{U}+\sum_{\beta+\gamma=\alpha}%
(n_{\beta,\gamma}c_{\beta}^{U}c_{\gamma}^{U}))ad(x_{\alpha}^{U})\text{,}%
\]
that is $\overline{\partial}_{0}\varphi_{\alpha}+\sum_{\beta+\gamma=\alpha
}(n_{\beta,\gamma}\varphi_{\beta} \wedge \varphi_{\gamma})=0$ for any $\alpha\in
\Phi^{+}$. Explicitly:%
\[
\left\{
\begin{tabular}
[c]{cc}%
$\overline{\partial}_{0}\varphi_{C_{i}}=0$ & $i=1,2\cdots,r$\\
$\overline{\partial}_{0}\varphi_{C_{i}+C_{j}}=n_{C_{i},C_{j}}\varphi_{C_{i}%
} \wedge \varphi_{C_{j}}$ & $\text{ if }C_{i}+C_{j}\in\Phi^{+}$\\
$\vdots$ &
\end{tabular}
\ \right.
\]

\bigskip
$\mathbf{Proposition  }$  Given any $(\varphi_{C_{i}})_{i=1}^{r}\in\Omega^{0,1}(Y, ~ \bigoplus_{i=1}^{n}\mathcal{O}(C_{i}))$ with $\overline{\partial}_{0}\varphi_{C_{i}}=0$
for every $i$, it can be extended to $\varphi=(\varphi_{\alpha})_{\alpha
\in\Phi^{+}}\in\Omega^{0,1}(Y,\bigoplus_{\alpha\in\Phi^{+}}\mathcal{O}(\alpha))$ such
that $\overline{\partial}_{\mathcal{\varphi}}^{2}=0$. Namely we have a
holomorphic vector bundle $\mathcal{E}_{\varphi}^{\mathfrak{g}}$ over $Y$.
\bigskip

This Proposition follows from the fact that for any $\alpha\in\Phi^{+}$, $H^{2}(Y,\mathcal{O}(\alpha))=0$. By computing the $A_n$, $D_n$, $E_6$, $E_7$ and $E_8$ types case by case, we have the following result:

\bigskip

$\mathbf{Theorem  }$ $\mathcal{E}_{\varphi}^{\mathfrak{g}}$ is trivial on $C_{i}$ if
and only if $[\varphi_{C_{i}}|_{C_{i}}]\neq0\in H^{1}(Y,\mathcal{O}_{C_{i}}(C_{i}))$.
\bigskip

The next lemma says that for any $C_{i}$, there always exists $\varphi
_{C_{i}}\in\Omega^{0,1}(Y,$ $\mathcal{O}(C_{i}))$ such that $0\neq\lbrack\varphi_{C_{i}%
}|_{C_{i}}]\in H^{1}(Y,$ $\mathcal{O}_{C_{i}}(C_{i}))\cong\mathbb{C}$.

\bigskip

$\mathbf{Lemma }$ For any $C_{i}$ in $Y$, the restriction homomorphism $H^{1}(Y,$
$\mathcal{O}_{Y}(C_{i}))\rightarrow H^{1}(Y,$ $\mathcal{O}_{C_{i}}(C_{i}))$ is surjective.

\section{Homogeneous $ADE$ bundles over flag varieties}
In this section, we will study the holomorphic structures of the homogeneous $ADE$ bundles $G \times \mathfrak{g}/B$ over $G/B$ when $\mathfrak{g}$ is of $ADE$ type.

\begin{lemma} For any finite dimensional representation $V$ of $B$, the associated representation bundle $G \times V /B$ over $G/B$ is an iterated extensions of holomorphic line bundles.

\begin{proof}
As $B$ is a solvable Lie group, using Lie's Theorem \cite{K}, any finite representation of $B$ has a filtration with irreducible factors. And any irreducible representation of $B$ is of one dimensional.
\end{proof}

\end{lemma}

Here we will first review the cohomology of line bundles over $G/B$, i.e. the Borel-Weil-Bott theorem \cite{Bott}\cite{D}. For the full flag variety $G/B$, we have $Pic(G/B)=\Lambda$, where $\Lambda$ is the weight lattice of the Lie algebra $\mathfrak{g}$. Hence for every $\lambda \in \Lambda$, we can associate a line bundle $L_{\lambda}$ over $G/B$. Denote

\[
\rho :=\frac{1}{2}\sum_{\alpha\in\Phi^{+}}\alpha=\sum_{i=1}^{r}\lambda_i
\]
where $\Phi^{+}$ is the set of positive roots of $\mathfrak{g}$ and $\{\lambda_1,\cdots,\lambda_r\}$ is the set of fundamental weights of $\mathfrak{g}$. Then by Borel-Weil-Bott theorem, we have

\bigskip

$\mathbf{(Borel-Weil-Bott  ~ Theorem)}$

(1) If $\lambda +\rho$ is singular (i.e. $\exists ~  \alpha \in \Phi$ such that ${\langle} {\alpha}^{\vee}, \lambda +\rho\rangle=0$), then
\[
H^i(G/B, L_{\lambda})=0 ~  \text{for all $i$};
\]

(2) If $\lambda +\rho$ is not singular, write $\lambda=\omega (\mu+\rho)-\rho$ with $\omega \in W$, $\mu \in C$, then
\[
H^i(G/B, L_{\lambda})=\left\{
\begin{tabular}
[c]{cc}%
$0$ & if $i\neq ind(\lambda +\rho)$\\
$V_{\mu}$ & if $i=ind(\lambda +\rho)$,\\
\end{tabular}
\ \right.
\]
where $\Phi$ is the set of roots of $\mathfrak{g}$, ${\alpha}^{\vee}$ is the dual root of $\alpha$, $W$ is the Weyl group, $C$ is the dominant chamber, $V_{\mu}$ is the irreducible representation of $G$ with highest weight $\mu$ and $H^i(G/B, L_{\lambda})$ is isomorphic to $V_{\mu}$ as $G$-modules, and for any $\lambda \in \Lambda$, $ind(\lambda)$ is defined to be the number of $\alpha\in\Phi^{+}$ such that $(\lambda, \alpha)<0$.

\bigskip

From the Borel-Weil-Bott theorem, we can compute some particular cases of cohomology of line bundles over $G/B$ easily.  Denote $\Delta=\{\alpha_1,\cdots,\alpha_r\}$ the set of simple roots.

\begin{proposition}
In $ADE$ cases, for any root $\alpha \in \Phi$, we have:

\[
(1) H^i(G/B, L_{\alpha})=0 ~  ~  ~  \text{for any} ~ i\geq 2;
\]

\[
(2) H^1(G/B, L_{\alpha})= \left\{
\begin{tabular}
[c]{cc}%
$\mathbb{C}$ & if $\alpha=-\alpha_i$ for some simple root $\alpha_i \in \Delta$\\
$0$ & otherwise.
\end{tabular}
\ \right .
\]

\begin{proof}
This proposition follows from the Borel-Weil-Bott theorem and the following lemma.
\end{proof}

\end{proposition}

\begin{lemma}
In $ADE$ cases, for any root $\alpha \in \Phi$, if $\alpha+\rho$ is not singular, then $ind(\alpha+\rho)\leq 1$.

\begin{proof}
Let $\Delta=\{\alpha_1,\cdots,\alpha_r\}$ be the set of simple roots and $\{\lambda_1,\cdots,\lambda_r\}$ be the set of corresponding fundamental weights. In $ADE$ cases, we have $(\alpha_i, \lambda_j)=\delta_{ij}$ for any $i, ~ j$ and for any roots $\alpha, ~ \beta \in \Phi$, we have $|(\alpha,\beta)|\leq2$ with $"="$ holds if and only if $\beta=\pm \alpha$.

If $\alpha \in \Phi$ with $\alpha \neq \pm \alpha_i$ for any $i$, then the coordinates of it in the basis of the fundamental weights are always $-1$, $0$ or $1$. Since $\rho :=\sum_{i=1}^{r}\lambda_i$ has all coordinates equal to $1$, $\alpha+\rho$ is either singular (a coordinate is $0$) or has index 0 (all coordinates positive).

If $\alpha= \alpha_i$ for some $i$, then the coordinates of it are always $-1$, $0$, $1$ or $2$. So $\alpha+\rho$ is either singular or has index 0.

If $\alpha=- \alpha_i$ for some $i$, for any $\beta \in \Phi^{+}$, if $(\alpha+\rho, \beta)<0$, then $\beta$ can only be $\alpha_i$, hence $ind(\alpha+\rho)= 1$.
\end{proof}
\end{lemma}

Now we will use this proposition to compute the holomorphic filtration structure of $G \times \mathfrak{g}/B$ when $\mathfrak{g}$ is of $ADE$ type.

\bigskip

\begin{example}
$\mathfrak{g}=A_n=sl(n+1, \mathbb{C})$, $G=SL(n+1, \mathbb{C})$ case. Choose $\Delta=\{x_1-x_2,\cdots,x_n-x_{n+1}\}$ as the set of simple roots. We first consider the associated representation bundle $G \times \mathbb{C}^{n+1} /B$. The representation $\mathbb{C}^{n+1}$ of $B$ has weights $\{x_1, \cdots, x_{n+1}\}$, with $\{v_1=(1,0,\cdots,0),\cdots,v_{n+1}=(0,\cdots,0,1)\}$ be their corresponding weight vectors. The filtration of this representation can only be
\[
\mathbb{C}^{n+1}\supset \mathbb{C}\langle v_2, v_3,\cdots, v_{n+1} \rangle \cdots \supset  \mathbb{C}\langle v_{n+1} \rangle \supset \{0\} .
\]
Hence the holomorphic structures of $G \times \mathbb{C}^{n+1} /B$ must be
\[
\overline{\partial}_{\varphi}=\left(
\begin{array}
[c]{cccc}%
\overline{\partial} & \varphi_{1,2} & \cdots & \varphi_{1,n+1}\\
0 & \overline{\partial} & \cdots & \varphi_{2,n+1}\\
\vdots & \vdots & \ddots & \vdots\\
0 & 0 & \cdots & \overline{\partial}%
\end{array}
\right)
\]
with $\varphi_{i,j}\in\Omega^{0,1}(G/B, L_{x_{n+2-i}}\otimes L_{x_{n+2-j}}^*)$ for any $j>i$. When $j>i$,
$x_{n+2-i}-x_{n+2-j}\in \Phi^-$ is a negative root, i.e. $\varphi_{i,j}\in\Omega^{0,1}(G/B, L_{\alpha})$ for some $\alpha \in \Phi^-$.

The integrability condition $\overline{\partial}_{\varphi}^{2}=0$ is
equivalent to, for $i=1,2,\cdots,n$,%
\[
\left\{
\begin{array}
[c]{l}%
\overline{\partial}\varphi_{i,i+1}=0,\\
\overline{\partial}\varphi_{i,j}=-\sum_{m=i+1}^{j-1}{\varphi}%
_{i,m} \wedge\varphi_{m,j}\text{, }~j\geq i+2\text{.}%
\end{array}
\right.
\]

Note $\varphi_{i,j}\in \Omega^{0,1}(G/B, L_{\alpha})$ for some $\alpha \in \Phi^-$. From%
\[
\sum_{m=i+1}^{j-1}[\varphi_{i,m}\wedge\varphi%
_{m,j}]\in H^{2}(G/B, L_{\alpha})=0\text{,}%
\]
we can find $\varphi_{i,j}$, such that $\overline{\partial}\varphi_{i,j}=-\sum_{m=i+1}^{j-1}{\varphi}
_{i,m} \wedge\varphi_{m,j}$.

Also $\overline{\partial}_{\varphi}^{2}=0$ tells us $\overline{\partial}\varphi_{i,i+1}=0$, i.e. $[\varphi_{i,i+1}] \in H^1(G/B, L_{x_{n+2-i}-x_{n+1-i}})\neq 0$ as $x_{n+1-i}-x_{n+2-i}$ is a simple root, hence we can take $[\varphi_{i,i+1}]$ to be a non-trivial class.

As $G \times \mathfrak{g}/B=G \times aut_0(\mathbb{C}^{n+1}) /B$, we have an induced holomorphic structure on $G \times \mathfrak{g}/B$ from $G \times \mathbb{C}^{n+1} /B$. From above, we can write the holomorphic structure of $G \times \mathbb{C}^{n+1} /B$ as $\overline{\partial}_{\varphi}:=\overline{\partial}_{0}+\sum_{\alpha\in\Phi^{-}}\rho(\varphi_{\alpha
})$, where $\rho$ is the representation $\mathfrak{g} \rightarrow End(\mathbb{C}^{n+1})$. More explicitly, if we write $\varphi
_{\alpha}=c_{\alpha}^{U}x_{\alpha}^{U}$ locally for some one form $c_{\alpha
}^{U}$ and the corresponding component $x_{\alpha}^{U}$ of locally Chevalley basis, then $\rho(\varphi_{\alpha})=c_{\alpha}^{U}\rho(x_{\alpha}^{U})$.
Now for a section $(x \in G/B, ~ X(x): \mathbb{C}^{n+1}\rightarrow \mathbb{C}^{n+1})$ of $G \times \mathfrak{g}/B$ and a section $(x \in G/B, ~ v(x) \in \mathbb{C}^{n+1})$ of $G \times \mathbb{C}^{n+1} /B$, we have

\begin{equation}
\nonumber
 \begin{aligned}[b]
\overline{\partial}_{\varphi}(X) \cdot v &= \overline{\partial}_{\varphi}(X\cdot v)-X \cdot (\overline{\partial}_{\varphi}v) \\
   &= (\overline{\partial}_{0}+\sum\rho(\varphi_{\alpha
}))(X \cdot v) -X \cdot (\overline{\partial}_{0}+\sum\rho(\varphi_{\alpha
}))\cdot v \\
   &= (\overline{\partial}_{0}X)\cdot v + \sum_{\alpha\in\Phi^{-}} c_{\alpha
} [\rho(x_{\alpha}), X]\cdot v \\
   &=(\overline{\partial}_{0}+\sum_{\alpha\in\Phi^{-}}c_{\alpha
} ad(x_{\alpha}))(X) \cdot v
 \end{aligned}
\end{equation}

Hence this induced holomorphic structure on $G \times \mathfrak{g} /B$ is
\[
\overline{\partial}_{\varphi}:=\overline{\partial}_{0}+\sum_{\alpha\in\Phi^{-}}ad(\varphi_{\alpha
})\text{,}%
\]
where $\overline{\partial}_{0}$ is the standard holomorphic structure and  $\varphi_{\alpha
} \in \Omega^{0,1}(G/B, L_{\alpha})$ for some $\alpha \in \Phi^-$.
\end{example}
\bigskip

\begin{example}
$\mathfrak{g}=D_n=so(2n, \mathbb{C})$, $G=SO(2n, \mathbb{C})$ case. Choose $\Delta=\{x_1-x_2, x_2-x_3,\cdots,x_{n-1}-x_{n}, x_{n-1}+x_{n}\}$ as the set of simple roots. We first consider the associated representation bundle $G \times \mathbb{C}^{2n} /B$. The representation $\mathbb{C}^{2n}$ of $B$ has weights $\{x_1, \cdots, x_n, x_{-n}, \cdots, x_{-1}\}$, with $\{v_1,\cdots,v_{2n}\}$ be their corresponding weight vectors. The filtration of this representation is not unique, in fact there are two choices, we arbitrary choose one:
\[
\mathbb{C}^{n+1}\supset \mathbb{C}\langle v_2, v_3,\cdots, v_{2n} \rangle \cdots \supset  \mathbb{C}\langle v_{2n} \rangle \supset \{0\} .
\]
Hence the holomorphic structures of $G \times \mathbb{C}^{2n} /B$ must be
\[
\overline{\partial}_{\varphi}=\left(
\begin{array}
[c]{cccc}%
\overline{\partial} & \varphi_{1,2} & \cdots & \varphi_{1,2n}\\
0 & \overline{\partial} & \cdots & \varphi_{2,2n}\\
\vdots & \vdots & \ddots & \vdots\\
0 & 0 & \cdots & \overline{\partial}%
\end{array}
\right)
\]
with $\varphi_{i,j}$ lies in $\Omega^{0,1}(G/B, L_{x_p-x_q})$ for some $p>q$, or $\Omega^{0,1}(G/B, L_{-x_p-x_q})$, or $\Omega^{0,1}(G/B, L_{-2x_p})$ when $i+j=2n+1$. Note that both $x_p-x_q$ and $-x_p-x_q$ are roots while $-2x_p$ is not.
Similar to the above $A_n$ case, to show that the integrability condition $\overline{\partial}_{\varphi}^{2}=0$ has solutions, we need to prove that $H^2(G/B, L_{-2x_p})=0$. The proof is similar to the proof of the above lemma, we will omit it here.

On the vector space $\mathbb{C}^{2n}$, we have a natural quadratic form $q$ such that $D_n=so(2n, \mathbb{C})=aut(\mathbb{C}^{2n},q)$, hence $G \times \mathfrak{g}/B=G \times aut(\mathbb{C}^{n+1}, q) /B$. To induce a holomorphic structure on $G \times \mathfrak{g}/B$ from $G \times \mathbb{C}^{2n} /B$, we need the holomorphic structure on  $G \times \mathbb{C}^{2n} /B$ to preserve the induced quadratic form $q$ on it. It is easy to check that $\overline{\partial}_{\varphi} q=0$ if and only if $\varphi_{i,j}=-\varphi_{2n+1-j,2n+1-i}$ for any $j>i$. From this, we know that all the nonzero $\varphi_{i,j}$'s are contained in $\Omega^{0,1}(G/B, L_{\alpha})$ for some $\alpha \in \Phi^-$. Furthermore, the induced holomorphic structure on $G \times \mathfrak{g} /B$ is
\[
\overline{\partial}_{\varphi}:=\overline{\partial}_{0}+\sum_{\alpha\in\Phi^{-}}ad(\varphi_{\alpha
})\text{,}%
\]
as in $A_n$ case. Note that these kind holomorphic structures don't depend on which filtration we choose for the representation at first.
\end{example}
\bigskip

Similar to the above examples, we can compute the holomorphic structures of $G \times \mathfrak{g}/B$ in $E_6$ and $E_7$ cases using the fact that $E_6=aut(\mathbb{C}^{27},c)$ and $E_7=aut(\mathbb{C}^{56},t)$, with some specific cubic form $c$ on $\mathbb{C}^{27}$ and quartic form $t$ on $\mathbb{C}^{56}$ (see \cite{CL} for more details). It turns out that, in these two cases, the induced holomorphic structures on $G \times \mathfrak{g}/B$ are also
\[
\overline{\partial}_{\varphi}:=\overline{\partial}_{0}+\sum_{\alpha\in\Phi^{-}}ad(\varphi_{\alpha
})\text{.}%
\]

Now we try to write the holomorphic structures on $G \times \mathfrak{g}/B$ directly. The filtration of the representation $\mathfrak{g}$ is given by the Chevalley order of its weights, hence not unique, we will choose an arbitrary one. Then the holomorphic structure on $G \times \mathfrak{g}/B$ can be written in a upper-triangular matrix as follows:
\[
\overline{\partial}_{\varphi}=\left(
\begin{array}
[c]{cccc}%
\overline{\partial} & \varphi_{1,2} & \cdots & \varphi_{1,N}\\
0 & \overline{\partial} & \cdots & \varphi_{2,N}\\
\vdots & \vdots & \ddots & \vdots\\
0 & 0 & \cdots & \overline{\partial}%
\end{array}
\right)
\]

\begin{proposition}
For the Lie algebra bundle $(G \times \mathfrak{g}/B, ~ [ ~ , ~ ])$ and $\overline{\partial}_{\varphi}$ as above, $\overline{\partial}_{\varphi} [ ~ , ~ ]=0$ if and only if $\overline{\partial}_{\varphi}=\overline{\partial}_{0}+\sum_{\alpha\in\Phi^{-}}ad(\varphi_{\alpha
})$ with  $\varphi_{\alpha
} \in \Omega^{0,1}(G/B, L_{\alpha})$ for some $\alpha \in \Phi^-$.

\begin{proof}
If $\overline{\partial}_{\varphi}=\overline{\partial}_{0}+\sum_{\alpha\in\Phi^{-}}ad(\varphi_{\alpha
})$, from Jacobi identity, we have $\overline{\partial}_{\varphi} [ ~ , ~ ]=0$.

Conversely, we suppose $\overline{\partial}_{\varphi} [ ~ , ~ ]=0$.

First, we can show that for those $\varphi_{i,j} \notin \Omega^{0,1}(G/B, L_{\alpha})$ for any $\alpha \in \Phi^-$, $\varphi_{i,j}=0$ by direct computations.

Second, for each $\alpha \in \Phi^-$, we consider those $\varphi_{i,j}$'s which are contained in $\Omega^{0,1}(G/B, L_{\alpha})$, it can be proved that these $\varphi_{i,j}$'s are different to each other by a constant coefficient.

Third, through more detailed calculations, we can see that for those $\varphi_{i,j}$'s contained in the same $\Omega^{0,1}(G/B, L_{\alpha})$, their impacts in $\overline{\partial}_{\varphi}$ is as $ad(\varphi_{\alpha
})$ for some $\varphi_{\alpha
} \in \Omega^{0,1}(G/B, L_{\alpha})$.
\end{proof}

\end{proposition}

From the above proposition, we know that the holomorphic structure $\overline{\partial}_{\varphi}$ of $G \times \mathfrak{g}/B$ do not depend on the filtration we choose at first. Now we consider the integrability condition of  $\overline{\partial}_{\varphi}=\overline{\partial}_{0}+\sum_{\alpha\in\Phi^{-}}ad(\varphi_{\alpha
})$. Since $H^{2}(G/B, L_{\alpha})=0$ for any $\alpha \in \Phi^-$, the integrability condition $\overline{\partial}_{\varphi}^{2}=0$ always has solutions, according to the computations in section $2$. Also $\overline{\partial}_{\varphi}^{2}=0$ tells us $\overline{\partial}\varphi_{-\alpha_i}=0$ for every simple root $\alpha_i$, i.e. $[\varphi_{-\alpha_i}] \in H^{1}(G/B, L_{-\alpha_i})\neq 0$, hence we can take $[\varphi_{-\alpha_i}]$ to be a non-trivial class. That means the holomorphic structure we got can be non-trivial.

As we mentioned, the associated bundle $G \times \mathfrak{g}/B$ is holomorphically trivial as the action of $B$ on $\mathfrak{g}$ can extend to the whole $G$. So the next question is what kind of $\overline{\partial}_{\varphi}$ can make $G \times \mathfrak{g}/B$ holomorphically trivial? To answer this question, we refer to the following theorem by X.Y. Pan in \cite{P}:

\bigskip

$\mathbf{Theorem  \cite{P}}$ For a homogenous space $G/P$, a vector bundle $V$ on $G/P$ is trivial if and only if the restriction of $V$ to every Schubert line is trivial.
\bigskip

Back to our cases, the Schubert lines in $G/B$ are given by $C_i=P_{\alpha_i}/B$, where $\alpha_i$ run through all the simple roots, and $P_{\alpha_i}$ is a parabolic subgroup of $G$ corresponding to $\alpha_i$.
\begin{lemma}
The bundle $(G \times \mathfrak{g}/B, \overline{\partial}_{\varphi}=\overline{\partial}_{0}+\sum_{\alpha\in\Phi^{-}}ad(\varphi_{\alpha
}))$ is holomorphically trivial if and only if $[\varphi_{-\alpha_i}|_{C_i}]\neq 0$ for every simple root $\alpha_i$.

\begin{proof}
Directly from \cite{CL}, $G \times \mathfrak{g}/B$ is trivial over $C_i=P_{\alpha_i}/B$ if and only if $[\varphi_{-\alpha_i}|_{C_i}]\neq 0$.
\end{proof}

\end{lemma}

The next lemma says that for any simple root $\alpha_{i}$, there always exists $\varphi
_{-\alpha_{i}}\in\Omega^{0,1}(G/B, L_{-\alpha_i})$ such that $\lbrack\varphi_{-\alpha_{i}%
}|_{C_{i}}]\in H^{1}(C_i, L_{-\alpha_i}|_{C_{i}})\cong H^{1}(\mathbb{P}^1, O(-2)) \cong \mathbb{C}$ is not zero.

\begin{lemma}
For any simple root $\alpha_{i}$, the restriction map $H^{1}(G/B, L_{-\alpha_i})\rightarrow H^{1}(C_i, L_{-\alpha_i}|_{C_{i}})$ is surjective.

\begin{proof}
From Borel-Weil-Bott theorem, we have $H^{1}(G/B, L_{-\alpha_i})\cong H^{0}(G/B, L_{0})$ as $S_{\alpha_i}(-\alpha_i+\rho)-\rho=0$. Also $H^{1}(C_i, L_{-\alpha_i}|_{C_{i}})\cong H^{0}(C_i, L_{0}|_{C_{i}})$ and the restriction map $H^{1}(G/B, L_{-\alpha_i})\rightarrow H^{1}(C_i, L_{-\alpha_i}|_{C_{i}})$ is the same with the restriction map $H^{0}(G/B, L_{0})\rightarrow H^{0}(C_i, L_{0}|_{C_{i}})$. From \cite{BS}\cite{Sen}, we know this restriction map is surjective.
\end{proof}
\end{lemma}

Since $H^{1}(G/B, L_{-\alpha_i})\cong H^{0}(G/B, L_{0}) \cong \mathbb{C}$, the above restriction map $H^{1}(G/B, L_{-\alpha_i})\rightarrow H^{1}(C_i, L_{-\alpha_i}|_{C_{i}})$ is in fact an isomorphism. Hence we have $\lbrack\varphi_{-\alpha_{i}
}|_{C_{i}}] \neq 0$ if and only if $\lbrack\varphi_{-\alpha_{i}
}] \neq 0$. Combine the above results, we have the following theorem:

\begin{theorem}
The holomorphic structure of $(G \times \mathfrak{g}/B, ~ [ ~ , ~ ])$ over $G/B$ is $\overline{\partial}_{\varphi}=\overline{\partial}_{0}+\sum_{\alpha\in\Phi^{-}}ad(\varphi_{\alpha
})$ with $[\varphi_{-\alpha_i}]\neq 0$ for every simple root $\alpha_{i}$.
\end{theorem}

Since $H^{1}(G/B, L_{-\alpha_i})\cong \mathbb{C}$ and $H^{1}(G/B, L_{-\alpha})=0$ for $\alpha \in \Phi^+$, $\alpha  \neq \alpha_i$, the holomorphic structure in the above theorem is unique up to isomorphism.

\section{$ADE$ bundles over cotangent bundles of the flag varieties}
In this section, we want to write the holomorphic structure of $G\times \mathfrak{n}\times \mathfrak{g}/B$ over $G\times \mathfrak{n}/B \cong T^*(G/B)$ when $\mathfrak{g}$ is of $ADE$ type. Similarly to Lemma 1, we know that $G\times \mathfrak{n}\times \mathfrak{g}/B$ is an iterated extensions of line bundles over $G\times \mathfrak{n}/B$ as $B$ is solvable. And any line bundle over $G\times \mathfrak{n}/B$ is the pull back of a line bundle over $G/B$ through the projection map $\pi: T^*(G/B)\cong G\times \mathfrak{n}/B \rightarrow G/B$. Denote $\mathfrak{L}_\lambda:=\pi^{*}L_\lambda$ to be the corresponding line bundle over $G\times \mathfrak{n}/B$ for any weight $\lambda \in \Lambda$.

Similar to the above section, we need to compute $H^1(G\times \mathfrak{n}/B, \mathfrak{L}_{\alpha})$ and $H^2(G\times \mathfrak{n}/B, \mathfrak{L}_{\alpha})$ for $\alpha \in \Phi^-$. We denote $H^i(\lambda):=H^i(G\times \mathfrak{n}/B, \mathfrak{L}_{\lambda})$ for convenience. Some properties and computations of  $H^i(\lambda)$ can be found in \cite{Bro}\cite{Bro2}\cite{H}.

Write $Cht(\lambda)$ for the combinatorial dimension of the interval $[\lambda^*, \lambda^+]$ in the Chevalley order (Here $\lambda^*$ is the unique dominant weight that is minimal with the property $\lambda^* \geq \lambda$ and $\lambda^+$ is the unique dominant weight on the Weyl group orbit of $\lambda$), i.e. the supremum over all $r$ such that there exists a chain
\[
\lambda^* \leq \mu_0 < \mu_1 < \cdots < \mu_r \leq \lambda^+
\]
with all $\mu_i$ dominant.

Various properties of $Cht(\lambda)$ can be found in \cite{Bro2}\cite{H}.
We recall the following from \cite{H} Lemma 4.2.

\bigskip

$\mathbf{Lemma  \cite{H}}$ (i) Let $\lambda \in \Lambda$, then $Cht(\lambda)=0$ iff $\lambda(\beta^\vee)\geq-1$ for all $\beta \in \Phi^+$. In particular, $Cht(\lambda)=0$ for all $\lambda \in C$.

(ii) Let $\lambda \in \Lambda$ with $Cht(\lambda)=0$ and let $\mu \in C$, then $Cht(\lambda+\mu)=0$.
\bigskip

We will also use the following theorem from \cite{Bro2}:

\bigskip
$\mathbf{Theorem  \cite{Bro2}}$ (i) For $\lambda \in \Lambda$, we have the equivalences
\[
H^i(\lambda)=0 ~ for ~ all ~ i\geq 1 \Leftrightarrow H^1(\lambda)=0 \Leftrightarrow Cht(\lambda)=0, i.e. ~ \lambda^*=\lambda^+.
\]

(ii) If $Cht(\lambda)=1$, then up to a shift in degrees
\[
H^1(\lambda)\simeq H^0(\lambda^*)/H^0(\lambda^+)[-ht(\lambda^+-\lambda^*)]\neq0.
\]

(iii) $H^i(\lambda)=0$ for $i>Cht(\lambda)$.
\bigskip

\begin{remark}
From the above lemma and theorem, in our $ADE$ cases, for any positive root $\alpha \in \Phi^+$, $Cht(\lambda)=0$, $H^i(\alpha)=0$ for all $i\geq 1$; for any negative root $\alpha \in \Phi^-$, $Cht(\lambda)\neq 0$, $H^1(\alpha)\neq 0$.
\end{remark}

\begin{proposition}
In our $ADE$ cases, for any negative root $\alpha \in \Phi^-$, $H^2(\alpha)=0$.
\end{proposition}

To prove this proposition, we need the following lemmas.

\bigskip
$\mathbf{Lemma \cite{Bro2}}$ Let $Q\supset P$ be two parabolic subgroups and let $V$ be an irreducible $P$-module. Write $Z:=G \times (\mathfrak{g}/\mathfrak{q})^*/P$.

(i) There exists at most one $i\geq 0$ such that
\[
H^i(Q/P, \mathfrak{L}_{Q/P}(V))\neq 0.
\]

(ii) If $H^i(Q/P, \mathfrak{L}_{Q/P}(V))=0$ for all $i\geq0$, then for all $i\geq0$,
\[
H^i(Z, \mathfrak{L}_{Z}(V))= 0.
\]

(iii) Suppose $\tilde{V}:=H^{v}(Q/P, \mathfrak{L}_{Q/P}(V))\neq 0$ for $v\geq0$, then
\[
H^i(Z, \mathfrak{L}_{Z}(V))=\left\{
\begin{tabular}
[c]{cc}%
$0$ & if $i<v$ \\
$H^{i-v}(Z, \mathfrak{L}_{Z}(\tilde{V}))$ & if $i\geq v$.
\end{tabular}
\ \right .
\]
Here $\mathfrak{L}_{Q/P}(V)$ and $\mathfrak{L}_{Z}(V)$ are the associated representation bundles over $Q/P$ and $Z$ respectively.
\bigskip

Now let $\alpha$ be a simple root, $Q=P_{\alpha} \supset B$, $X=T^*(G/B)=G\times \mathfrak{n}/B$. Then $\mathfrak{L}_{-\alpha}=\mathfrak{L}_{X}(\mathbb{C}_{\alpha})^*$ on $X$ has a natural linear section with scheme of zeros $Z=G \times (\mathfrak{g}/\mathfrak{q})^*/B$, hence we can identify $\mathfrak{L}_{X}(\mathbb{C}_{\alpha})[-1]$ with the ideal sheaf of $Z$ in $X$, where the $[-1]$ denotes a shift in grading such that generators have degree $1$. Write $\iota: Z\subset X$ for the inclusion, then for any weight $\lambda$, we have a $G$-equivariant  exact sequence of graded $\mathcal{O}_{X}$-modules
\[
0 \rightarrow \mathfrak{L}_{X}(\mathbb{C}_{\lambda+\alpha})[-1] \rightarrow \mathfrak{L}_{X}(\mathbb{C}_{\lambda}) \rightarrow \iota_*\mathfrak{L}_{Z}(\mathbb{C}_{\lambda})\rightarrow 0.
\]

As before, we write $H^i(\lambda):=H^i(X, \mathfrak{L}_{X}(\mathbb{C}_{\lambda}))$ and $H_{\alpha}^i(\lambda):=H^i(Z, \mathfrak{L}_{Z}(\mathbb{C}_{\lambda}))$, then we have a long exact sequence
\[
\cdots \rightarrow H^i(\lambda+\alpha)[-1] \rightarrow H^i(\lambda) \rightarrow H_{\alpha}^i(\lambda) \rightarrow \cdots
\]

\begin{lemma}
For any simple root $\alpha$ and any weight $\lambda$, if $\langle \lambda, \alpha \rangle =-1$, then $H^i(\lambda) \cong H^i(\lambda+\alpha)[-1]$ for any $i\geq0$.

\begin{proof}
Write $Q=P_{\alpha} \supset B$, then $H^i(Q/B, \mathfrak{L}_{Q/B}(\mathbb{C}_{\lambda})) \cong H^i(\mathbb{P}^1, \mathcal{O}(-1))=0$ for any $i\geq0$. From (ii) of the above Lemma \cite{Bro2}, we have $H_{\alpha}^i(\lambda)=0$ for any $i\geq0$. Hence $H^i(\lambda) \cong H^i(\lambda+\alpha)[-1]$ for any $i\geq0$ by the above long exact sequence.
\end{proof}
\end{lemma}

\begin{lemma} For any simple root $\alpha$, $H^2(-\alpha)=0$.

\begin{proof}
Consider the above long exact sequence
\[
\cdots \rightarrow H^i(\lambda+\alpha)[-1] \rightarrow H^i(\lambda) \rightarrow H_{\alpha}^i(\lambda) \rightarrow \cdots
\]

Take $\lambda=-\alpha$, since $H^2(0)=0$ (directly from Theorem  \cite{Bro2} (i) and $Cht(0)=0$), to show $H^2(-\alpha)=0$, we only need to show $H_{\alpha}^2(-\alpha)=0$.

Take $\lambda=0$, since $H^1(0)=0$ and $H^2(\alpha)=0$, we have $H_{\alpha}^1(0)=0$.

From Lemma \cite{Bro2}, since $\tilde{V}:=H^{1}(P_{\alpha}/B, \mathfrak{L}_{P_{\alpha}/B}(\mathbb{C}_{-\alpha}))\cong H^1(\mathbb{P}^1, \mathcal{O}(-2)) \cong \mathbb{C} \neq 0$ and the action of $B$ on $\tilde{V}$ is trivial (from Borel-Weil-Bott theorem), we have $H_{\alpha}^2(-\alpha)=H_{\alpha}^1(0)=0$.
\end{proof}
\end{lemma}

From the above two lemmas, we can prove our Proposition 11 now.

\begin{proof} ($\mathbf{Proposition ~ 11}$) We want to prove $H^2(\lambda)=0$ for any negative root $\lambda \in \Phi^-$ in our $ADE$ cases.

If $ht(\lambda)=1$, i.e. $\lambda=-\alpha$ for some simple root $\alpha$, then from Lemma 13, $H^2(\lambda)=0$.

By induction on $ht(\lambda)$. Suppose the proposition is true for every $\beta \in \Phi^-$ with $ht(\beta)=m$. Given any $\lambda \in \Phi^-$ with $ht(\lambda)=m+1$, by Lemma A in section 10.2 of \cite{H2}, there exists some simple root $\alpha$ such that $\langle \lambda, \alpha \rangle =-1$, i.e. $\lambda+\alpha \in \Phi^-$ with $ht(\lambda+\alpha)=m$, hence $H^2(\lambda+\alpha)=0$. From Lemma 12, we have $H^2(\lambda)=H^2(\lambda+\alpha)=0$.
\end{proof}

\bigskip

As in the above section, we now try to write the holomorphic structures $\overline{\partial}_{\varphi}$ on $G\times \mathfrak{n}\times \mathfrak{g}/B$ directly. Since $G\times \mathfrak{n}\times \mathfrak{g}/B$ is an iterated extensions of line bundles, $\overline{\partial}_{\varphi}$ can be written in a upper-triangular matrix, depending on the filtrations we choose for the representation $\mathfrak{g}$. For $\overline{\partial}_{\varphi}$ to preserve the Lie bracket, $\overline{\partial}_{\varphi}$ can only be $\overline{\partial}_{\varphi}=\overline{\partial}_{0}+\sum_{\alpha\in\Phi^{-}}ad(\varphi_{\alpha
})$ with $\varphi_{\alpha} \in \Omega^{0,1}(G\times \mathfrak{n}/B, \mathfrak{L}_{\alpha})$ for some $\alpha \in \Phi^-$. Hence $\overline{\partial}_{\varphi}$'s are not depending on the filtrations we choose at first. For the integrability condition $\overline{\partial}_{\varphi}^{2}=0$ to have solutions, we need $H^2(\alpha)=0$ for any negative root $\alpha \in \Phi^-$, which is true by Proposition 11. Also $\overline{\partial}_{\varphi}^{2}=0$ tells us $\overline{\partial}\varphi_{-\alpha_i}=0$ for every simple root $\alpha_i$, i.e. $[\varphi_{-\alpha_i}] \in H^{1}(G\times \mathfrak{n}/B, \mathfrak{L}_{-\alpha_i})\neq 0$, hence we can take $[\varphi_{-\alpha_i}]$ to be a non-trivial class. That means the holomorphic structures we got can be non-trivial.

Similar to $G\times \mathfrak{g}/B$ case, the associated bundle $G\times \mathfrak{n}\times \mathfrak{g}/B$ is holomorphically trivial as the action of $B$ on $\mathfrak{g}$ can extend to the whole $G$. For example, we can take the holomorphic structure of $G\times \mathfrak{n}\times \mathfrak{g}/B$ to be $\pi^{*}(\overline{\partial}_{\varphi})$ where $\pi : G\times \mathfrak{n}\times \mathfrak{g}/B \rightarrow G\times \mathfrak{g}/B$ is the projection map and $\overline{\partial}_{\varphi}$ is the holomorphic structure of $G\times \mathfrak{g}/B$ as in Theorem 9. That means $G\times \mathfrak{n}\times \mathfrak{g}/B$ is a pull back of $G\times \mathfrak{g}/B$ holomorphically, hence trivial. In general, for $G\times \mathfrak{n}\times \mathfrak{g}/B$ over $G\times \mathfrak{n}/B$ to be trivial, its restriction to $G/B$ must also be trivial, hence for each simple root $\alpha_i$, $[\varphi_{-\alpha_i}|_{G/B}] \neq 0$. As $H^{1}(G\times \mathfrak{n}/B, \mathfrak{L}_{-\alpha_i})=\bigoplus_{j=0}^{\infty}H^{1}(G/B,S^{j}\mathfrak{n}^{*}\otimes L_{-\alpha_i})$, where $S^{j}\mathfrak{n}^{*}=G \times s^{j}\mathfrak{n}^{*}/B$ is the associated vector bundle over $G/B$ and $s^{j}\mathfrak{n}^{*}$ is the $j$-th symmetric power of the dual space $\mathfrak{n}^{*}$ of $\mathfrak{n}$, the restriction map $H^{1}(G\times \mathfrak{n}/B, \mathfrak{L}_{-\alpha_i}) \rightarrow H^{1}(G/B, L_{-\alpha_i})$ is just the projection, hence it is surjective. That means we can always take $[\varphi_{-\alpha_i}] \in H^{1}(G\times \mathfrak{n}/B, \mathfrak{L}_{-\alpha_i})$ such that $[\varphi_{-\alpha_i}|_{G/B}] \neq 0$.

Combine the above results, we have the following theorem:

\begin{theorem}
The holomorphic structure of $(G \times \mathfrak{n} \times \mathfrak{g}/B, ~ [ ~ , ~ ])$ over $G \times \mathfrak{n}/B$ is $\overline{\partial}_{\varphi}=\overline{\partial}_{0}+\sum_{\alpha\in\Phi^{-}}ad(\varphi_{\alpha
})$ with $[\varphi_{-\alpha_i}|_{G/B}] \neq 0$ for every simple root $\alpha_i$.
\end{theorem}

\section{$ADE$ bundles over $ADE$ singular surfaces}
In this section, we consider the restriction of the $\mathfrak{g}$-bundle $G \times \mathfrak{n} \times \mathfrak{g}/B$ from $G \times \mathfrak{n} /B$ to $\widetilde{S}$, note that $\widetilde{S}$ is the minimal resolution of the $ADE$ singular surface $S=\mathbb{C}^2/\Gamma$. It is obviously that this $\mathfrak{g}$-bundle over $\widetilde{S}$ is also an iterated extensions of line bundles.

Denote $C=\bigcup C_i$ to be the exceptional locus in $\widetilde{S}$, with each $C_i$ irreducible component, then the dual graph of $C$ is an $ADE$ Dynkin diagram of the corresponding type. The Picard group of $\widetilde{S}$ is a free abelian group generated by divisors dual to the irreducible curves $C_i$ \cite{M}, i.e. $Pic(\widetilde{S})=\mathbb{Z}\langle D_i\rangle$ with each $D_i$ dual to $C_i$.

As before, we know that the irreducible curves $C_i=P_{\alpha_i}/B$ are Schubert lines in $G/B$, where $\alpha_i$ run through all the simple roots. Now for any weight $\lambda$, we consider the restriction of the line bundle $L_{\lambda}$ from $G/B$ to $C_i$, it is easy to see that $L_{\lambda}|_{C_i}\cong \mathcal{O}_{\mathbb{P}^1}(\langle \lambda, \alpha_i \rangle)$. How about the restriction of the line bundle $\mathfrak{L}_\lambda=\pi^{*}L_\lambda$ form $G \times \mathfrak{n} /B$ to $\widetilde{S}$?

\begin{lemma}
For any root $\alpha=\sum n_i\alpha_i$, $\mathfrak{L}_{\alpha}|_{\widetilde{S}}\cong\mathcal{O}_{\widetilde{S}}(\sum -n_iC_i)$.

\begin{proof}
For the simplicity of computations, we first assume $\mathfrak{L}_{\alpha}|_{\widetilde{S}}\cong\mathcal{O}_{\widetilde{S}}(\sum -m_iC_i)$ with $m_i$'s integers.

We consider $\mathfrak{L}_{\alpha}|_{C_j}$ for each $j$, then
\[
(-\sum m_iC_i) \cdot C_j=\langle \sum n_i\alpha_i, \alpha_j \rangle
\]
i.e.
\[
\sum m_i(C_i \cdot C_j)=-\sum n_i \langle \alpha_i, \alpha_j \rangle
\]
Since $[C_i \cdot C_j]_{r\times r}=[-\langle \alpha_i, \alpha_j \rangle]_{r\times r}$ are invertible matrices, here $r$ is the rank of the Lie algebra $\mathfrak{g}$,
\[
(m_1,\cdots,m_r)[C_i \cdot C_j]=(n_1,\cdots,n_r)[-\langle \alpha_i, \alpha_j \rangle]
\]
has unique solution $(m_1,\cdots,m_r)=(n_1,\cdots,n_r)$. Hence our assumption is right and $\mathfrak{L}_{\alpha}|_{\widetilde{S}}\cong\mathcal{O}_{\widetilde{S}}(\sum -n_iC_i)$.
\end{proof}
\end{lemma}

From the above lemma, the $\mathfrak{g}$-bundle over $\widetilde{S}$ topologically is
\[
\mathcal{O}^{\oplus r} \oplus \bigoplus_{(\sum n_i C_i)^2=-2}\mathcal{O}(\sum n_i C_i)
\]

Since the holomorphic structure over $G \times \mathfrak{n} \times \mathfrak{g}/B$ is $\overline{\partial}_{\varphi}=\overline{\partial}_{0}+\sum_{\alpha\in\Phi^{-}}ad(\varphi_{\alpha
})$, the induced holomorphic structure of this $\mathfrak{g}$-bundle over $\widetilde{S}$ is also
\[
\overline{\partial}_{\varphi}=\overline{\partial}_{0}+\sum_{\alpha\in\Phi^{-}}ad(\varphi_{\alpha
})
\]
where for each $\alpha=\sum -n_i\alpha_i \in \Phi^-$, $\varphi_{\alpha} \in \Omega^{0,1}(\widetilde{S}, \mathcal{O}_{\widetilde{S}}(\sum n_iC_i))$.

From the rationality of $\widetilde{S}$, we have $H^1(\widetilde{S}, \mathcal{O})=H^2(\widetilde{S}, \mathcal{O})=0$, hence $H^1(\widetilde{S}, \mathcal{O}(C_i)) \cong \mathbb{C}$, $H^2(\widetilde{S}, \mathcal{O}(C_i))=0$ and the restriction map $H^1(\widetilde{S}, \mathcal{O}(C_i)) \rightarrow H^1(\widetilde{S}, \mathcal{O}_{C_i}(C_i))\cong \mathbb{C}$ is an isomorphism, for every $C_i$. Similar to the proof of Proposition 11, using Lemma A in section 10.2 of \cite{H2} and induction on $ht(\sum n_iC_i)$, we can show that for each effective divisor $D=\sum n_iC_i$ with $D^2=-2$, $H^2(\widetilde{S}, \mathcal{O}(D))=0$. This implies that the integrability condition $\overline{\partial}_{\varphi}^{2}=0$ always has solutions. Also $\overline{\partial}_{\varphi}^{2}=0$ tells us $\overline{\partial}\varphi_{-\alpha_i}=0$ for every simple root $\alpha_i$, i.e. $[\varphi_{-\alpha_i}] \in H^{1}(\widetilde{S}, \mathcal{O}(C_i))\neq 0$, hence we can take $[\varphi_{-\alpha_i}]$ to be a non-trivial class. That means the holomorphic structures we got can be non-trivial. For the $\mathfrak{g}$-bundle to be trivial over $\widetilde{S}$, it must be trivial over each $C_i$, hence $[\varphi_{-\alpha_i}|_{C_i}] \neq 0$, which is the same with $[\varphi_{-\alpha_i}] \neq 0$.

\begin{theorem}
The restriction of the $\mathfrak{g}$-bundle $G \times \mathfrak{n} \times \mathfrak{g}/B$ from $G \times \mathfrak{n} /B$ to $\widetilde{S}$ is
 \[
(\mathcal{O}^{\oplus r} \oplus \bigoplus_{(\sum n_i C_i)^2=-2}\mathcal{O}(\sum n_i C_i),  ~ \overline{\partial}_{\varphi}=\overline{\partial}_{0}+\sum_{\alpha\in\Phi^{-}}ad(\varphi_{\alpha
})).
\]
with $[\varphi_{-\alpha_i}] \neq 0$ for every simple root $\alpha_i$.
\end{theorem}

\bigskip

We can easily note that the holomorphic structures here have the same form with the holomorphic structures in \cite{CL}.

\bigskip

School of Science, East China University of Science and Technology, Meilong
Road 130, Shanghai, China

E-mail address: yxchen76@ecust.edu.cn

\bigskip

The Institute of Mathematical Sciences and Department of Mathematics, The
Chinese University of Hong Kong, Shatin, N.T., Hong Kong

E-mail address: leung@math.cuhk.edu.hk

\end{document}